\documentclass[review]{elsarticle}
\usepackage{lineno,hyperref}

\usepackage{graphicx}
\usepackage[all]{xy}
\usepackage{bm}
\usepackage{amscd,amssymb,amsbsy,amsmath,amsfonts}
\usepackage{color}
\usepackage{fancybox}

 \usepackage{graphics}%,epstopdf}
\newtheorem{theorem}{Theorem}

\newtheorem{definition}{Definition}

\modulolinenumbers[5]

\journal{MATUA}

\begin{document}


\begin{frontmatter}

%\dochead{xxx-xxxx}
%\titulo{El t\'itulo en espa\~ol}
%\titulo{De Poincar\'e a May: El nacimiento de la Din\'amica Discreta}
\title{From Poincar\'e to May: The Genesis of Discrete Dynamics}


\author{Oscar Eduardo Mart\'inez Castiblanco}

\address{Universidad Sergio Arboleda}
\ead{oscar.martinez@usa.edu.co}

\author{Primitivo Acosta-Hum\'anez}

\address{Instituto Superior de Formaci\'on Docente Salom\'e Ure\~{n}a - ISFODOSU}
\ead{primitivo.acosta-humanez@isfodosu.edu.do}

\begin{abstract}
In this paper, the origin of discrete dynamics is stated from a historical point of view, as well as its main ideas: fixed and periodic points, chaotic behaviour, bifurcations. This travel will begin with Poincar\'e's work and will finish with May's work, one of the most important scientific papers of 20th century.

This paper is based on the M.Sc. thesis "Simple Permutations, Pasting and Reversing" (\cite{Tesis}), written by the first author under the guidance of the second author.
\end{abstract}

\end{frontmatter}
\vspace*{10pt}
%\begin{multicols}{2}

% \linenumbers


\section{Introduction}

This paper corresponds to the theoretical framework of the master thesis \cite{Tesis}, in where were combined techiques of dynamical systems (some of them will be presented here) with \emph{Pasting and Reversing} techiques (see \cite{Genealogy,Pegamiento,PastingMatrix,preprint,Observaciones,Simple} ). In particular, in this paper the origin of Discrete Dynamics will be presented. A meaningful context may be useful to understand the relevance and further applications of this branch of Dynamical Systems. In this way, a historical landscape will help the reader to answer the question: Why and what for do we study discrete dynamical systems?

\section{Poincar\'e}

Going back in time, the origin of this problem can be set in one name: Henri Poincar\'e. Born into a middle-upper class family, Poincar\'e was one of the last universal mathematicians. He entered the \'Ecole Polytechnique in 1873, receiving his doctorate in 1879, after this, began a university career at University of Caen and then, in 1881, at University of Paris (\cite{History}). Poincar\'e's work includes contributions to almost all areas of mathematics, even physics and theoretical astronomy.

In fact, The beginning of this story lies on an astronomical problem: The three-body problem. It is well known that a system of two orbiting masses, interacting through gravitational acceleration can be described by using a differential equation. In this case, the system is analytically solvable. However, a system with three or more masses under the same hypotheses is not analytically solvable.

Since there is no way to determine the significant number of state spaces analytically, A contest was held to produce the best research in celestial mechanics, related to the stability of the solar system (a particular case of the n-body problem). This contest, held in 1889 to commemorate the 60th birthday of King Oscar II of Sweden and Norway (\cite{Chaos}) declared Poincar\'e as winner. 

In a field mainly dominated by quantitative methods so far (i.e. series expansions), Poincar\'e innovated by using new quantitative methods and simplifying assumptions. He assumed that the three bodies can be in a plane, since this decreases one degree of freedom. Furthermore, he assumed that one of the masses can be depreciated respect to the other two (e.g. Earth, Moon, and an artificial satellite orbiting among them). From this point of view, the problem turns into describing the third (small) body orbit using as reference the other two. Even when it was well known that two large masses would travel in ellipses, Poincar\'e considered once again a particular case: they would move in circles (with center in the center of mass of the system) at a constant speed. This simpler way to face the problem,took him into the ideas of stable and unstable manifolds, as well as homoclinic points.

\subsection{Poincar\'e Sections}

Even if all these ideas are revolutionary at that time, there is still one more, an idea that will create a whole new way to understand dynamic systems: The Poincar\'e section. This method is used to look complex trajectories at a simpler way. Instead of study (trace) the whole trajectory, Poincar\'e consider the intersection of it with a two-dimensional plane. This reduces the problem from $n$ to $n-1$ dimensions. Besides, in this particular case, the dynamic system is reduced into a discrete map.

\begin{definition}
A map is a formula that describes the new state in terms of the previous state.
\end{definition}

Consider $C$ as a trajectory (i.e. the solution of a differential equation with initial conditions A) and $S$ as a Two-dimensional plane. The intersections of $C$ and $S$ determines a (discrete) set of points. This set is the Poincar\'e map.
\\
\begin{definition}
Let G be A map defined by $C \cap S$, such that $G(A)=B$. G is a Poincar\'e map.
\end{definition}

\begin{figure}
\centering
\includegraphics[scale=0.25]{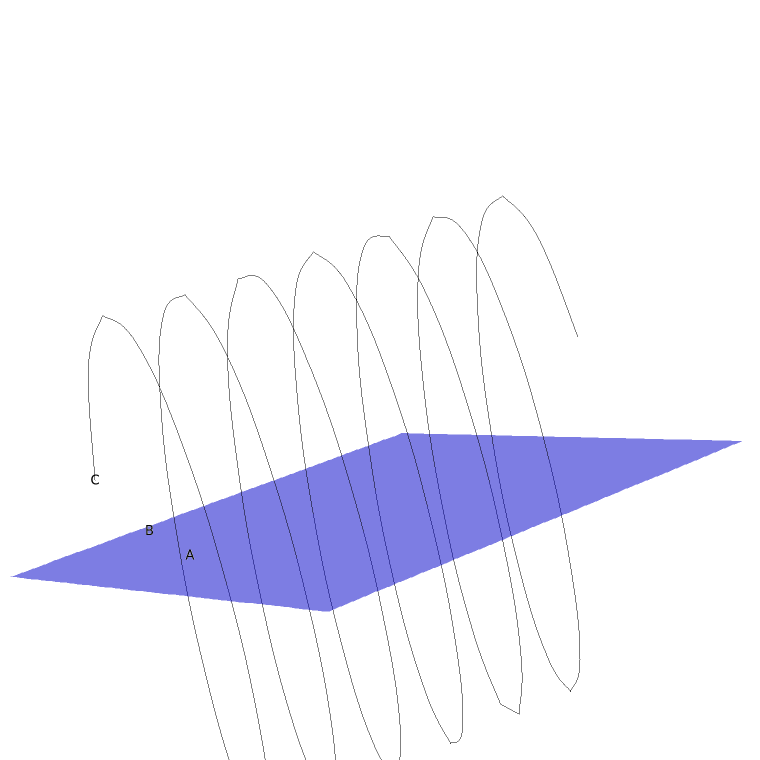}
\caption[Original Poincar\'e Section]{Poincar\'e section}
\label{Poincare2}
\end{figure}

The plane $S$ is named surface of section. By using this technique, Chaotic  behavior of differential equations can be studied by "reduction" to discrete dynamics. An orbit will be periodic if there exists $n$ such that $G^{(n)}(A)=A$.

\begin{figure}
\centering
\includegraphics[scale=0.75]{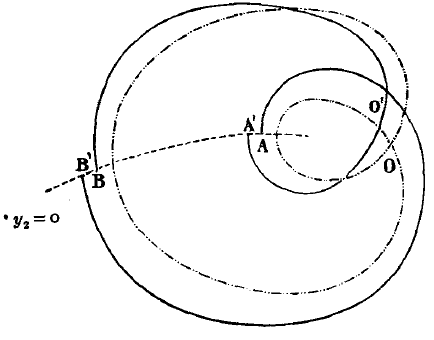}
\caption[Original Poincar\'e Section]{Poincar\'e section as in \cite{Poincare1}}
\label{Poincare1}
\end{figure}
Poincar\'e sections were applied in problems related with the \emph{integrability and non-integrability of hamiltonian dynamical systems}, see \cite{emalca,siam,AABD,AAD,AAS}
Since Poincar\'e maps are simpler to evaluate, even by computer, they were used to study complex trajectories, such as the ones determined by Lorenz system in (\cite{Lorenz}). Once the discrete dynamics is created, our path follows in two branches: Lorenz (and H\'enon) and May.

\section{Lorenz and H\'enon}

Edward Lorenz (1917 - 2008) began his career as mathematician, however, during WWII he moved into Atmospheric Sciences (Weather predictions) (\cite{Blanchard}). He studied systems with forced dissipative hydrodynamic flow, whose solutions can be understood as trajectories in phase space, and based on these solutions, examine the feasibility of very-long-range weather predictions. If it is true that some hydrodynamical systems show either steady-state flow patterns or periodic oscillations, Lorenz found a case (cellular convection) in which all solutions are unstable and nonperiodic.

\begin{figure}
\centering
\includegraphics[scale=1]{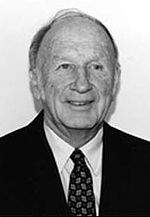}
\caption[Lorenz portrait]{Edward Lorenz as in \cite{Blanchard}}
\label{Poincare1}
\end{figure}

\subsection{Lorenz Chaotic Attractor}

In Deterministic Nonperiodic Flow (\cite{Lorenz}) Lorenz stated instability conditions with respect to modifications of small amplitude. This instability is due to sensitive dependence on initial conditions. 
Let $P(t)$ a trajectory. The trajectory is unstable if $\mid P(t_{1})-P(t_{1} + \tau)\mid < \varepsilon$ for some $t$ but $\mid P(t)-P(t + \tau)\mid \nless \varepsilon$ as $t \rightarrow \infty$. 

\begin{definition}
The Lorenz system 
\[ \begin{cases} 
      \dot{X}= & - \sigma X + \sigma Y \\
      \dot{Y}= & -XZ + rX - Y \\
      \dot{Z}= & XY - bZ 
   \end{cases}
\]
\end{definition}

The stability of the system can be understood by linearization

\[
\left[ 
\begin{array}{c}
x_0\\
y_0\\
z_0\\
\end{array}
\right]
=
\begin{bmatrix}
- \sigma & \sigma & 0\\
(r-Z) & -1 & -X\\
Y & X & -b
\end{bmatrix}
\left[ 
\begin{array}{c}
x_0\\
y_0\\
z_0\\
\end{array}
\right] 
\]

\begin{figure}
\centering
\includegraphics[scale=0.25]{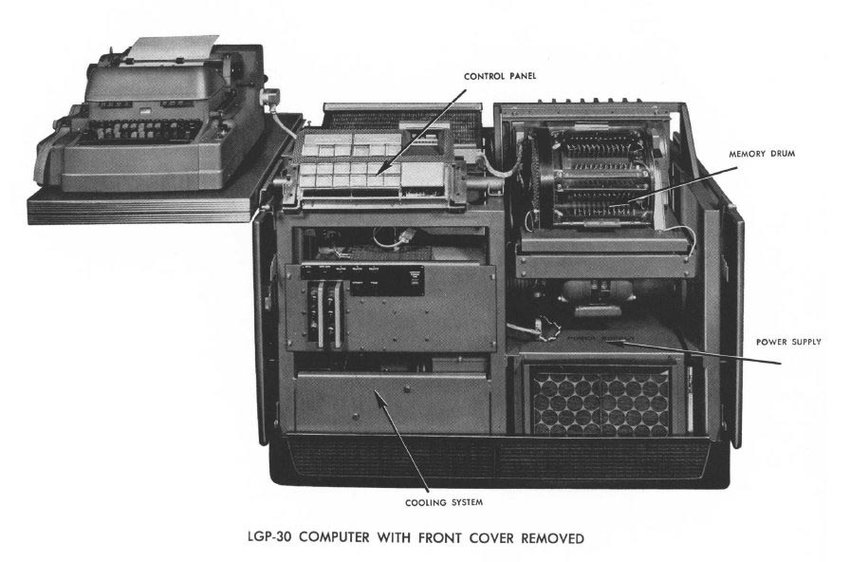}
\caption{Royal McBee LGP-30 electronic computer}
\label{Royal}
\end{figure}

This system has a steady-state solution if $X=Y=Z=0$. With this solution, the characteristic equation of this matrix is 
\[
[\lambda + b][\lambda^2 + (\sigma + 1) \lambda + \sigma (1-r)]=0 
\]
This equation has three real roots if $r>0$, and one of them is positive if $r>1$. In this case the system has two more steady-state solutions $X=Y=\pm \sqrt{b(r-1)},Z=r-1$. Taking these solutions into the linearization system the new characteristic equation is 
\[
\lambda^3 + (\sigma+b+1)\lambda^2 + (r+\sigma)b\lambda + 2\sigma b(r-1)=0
\]
Once again, if $r>1$ this equation has one real negative root. Besides, if 
\[
r=\frac{\sigma(\sigma + b + 3)}{\sigma - b - 1}
\]
 the complex conjugate roots are pure imaginary. Unfortunately, this information about Lorenz system obtained by linearization only applies over small perturbations, so numerical integration is required. Following Saltzman's work, Lorenz used $\sigma = 10, r=28, b=\frac{8}{3}$, a double-approximation procedure, and a Royal McBee LGP-30 electronic computer  to obtain numerical solutions (see Figure \ref{Royal}).

In this method, an initial time $t_0$ and an increment $\Delta t$ are set, and let $X_{i,n} = X_{i}(t_0+n \Delta t)$. Then we introduce the auxiliary approximations
\begin{eqnarray*}
X_{i,n+1} = X_{i,n} + F_i (P_n) \Delta t \\
X_{i,n+2} = X_{i,n+1} + F_i (P_{n+1}) \Delta t \\
\end{eqnarray*}
where $P_n = (X_{1,n},\ldots,X_{M,n})$, with $M$ as the number of variables of the system. The double-approximation procedure is defined by 
\[
X_{i,n+1}=X_{i,n}+\frac{1}{2}[F_i (P_n) + F_i (P_{n+1})] \Delta t
\]

From the previous equations it follows that the double-approximation procedure can be rewritten as:
\begin{eqnarray*}
X_{i,n+1} & = & X_{i,n}+\frac{1}{2}[F_i (P_n) + F_i (P_{n+1})] \Delta t \\
& = &X_{i,n}+\frac{1}{2}[\frac{X_{i,n+1}-X_{i,n}}{\Delta t} + \frac{X_{i,n+2}-X_{i,n+1}}{\Delta t}] \Delta t \\
& = &X_{i,n}+\frac{1}{2}[X_{i,n+1}-X_{i,n} + X_{i,n+2} -X_{i,n+1}\\
& = & \frac{1}{2}[X_{i,n+2} + X_{i,n}]\\
\end{eqnarray*}

Lorenz concluded that "when our results concerning the instability of nonperiodic flow are applied to the atmosphere, which is ostensibly nonperiodic, they indicate that prediction of the sufficiently distant future is impossible by any method" (\cite{Lorenz}). However, this is not the key point in our story. He found that all solutions are confined within the same bounds. This object that "traps" all solutions is known as Lorenz Attractor.

\begin{figure}
\centering
\begin{tabular}{c c c}
\includegraphics[scale=0.4]{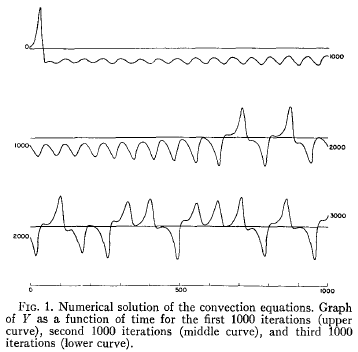}  &
\includegraphics[scale=0.4]{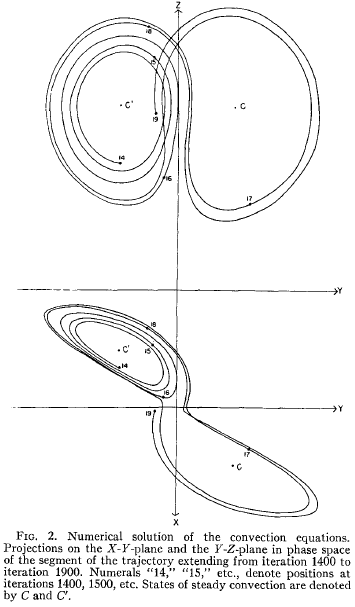}  & \includegraphics[scale=0.4]{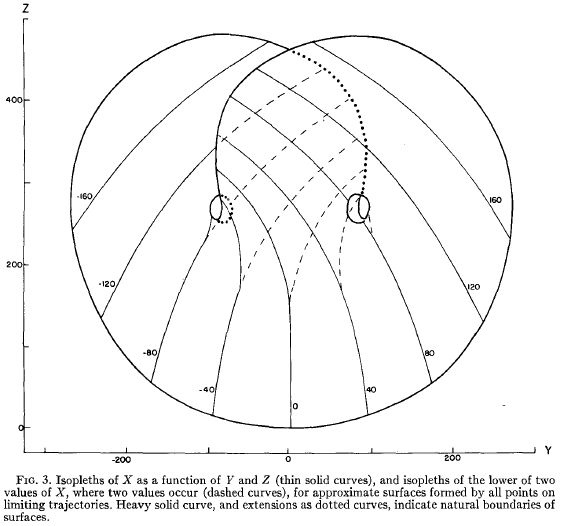}\\
\end{tabular}
\caption[AtractorOriginal]{Numerical solutions of Lorenz System, as in \cite{Lorenz}}
\end{figure}

As we can see in Lorenz paper, he could catch the idea of the "butterfly" trapping all trajectories. In fact, he did a remarkable work showing the geometry of the attractor with the tools at his reach, by using phase portraits and particular trajectories.

\begin{figure}
\centering
\includegraphics[scale=0.5]{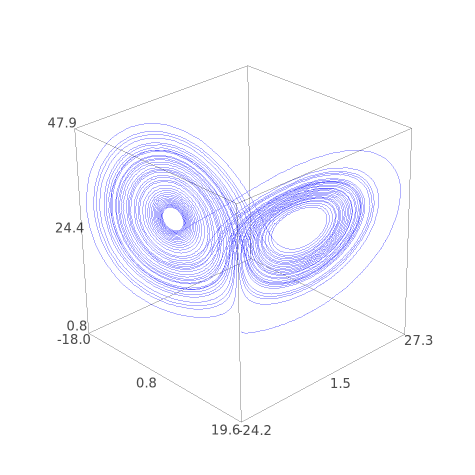}
\caption[AtractorNow]{Lorenz attractor for $\sigma = 10, r=28, b=\frac{8}{3}$}
\label{AttractorNow}
\end{figure}

The attractor (as it can be seen in 3D plots) is a complex object, and due to its properties (explained previously), the use or Poincar\'e sections to study it seems to be a good idea. Besides, the use of computers favour this approach to the problem.

\begin{figure}
\centering
\begin{tabular}{c c}
\includegraphics[scale=0.25]{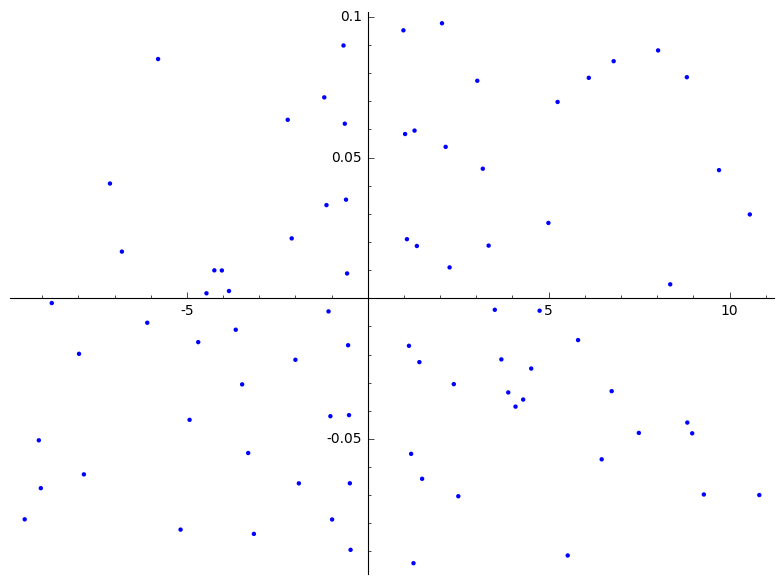}  & \includegraphics[scale=0.25]{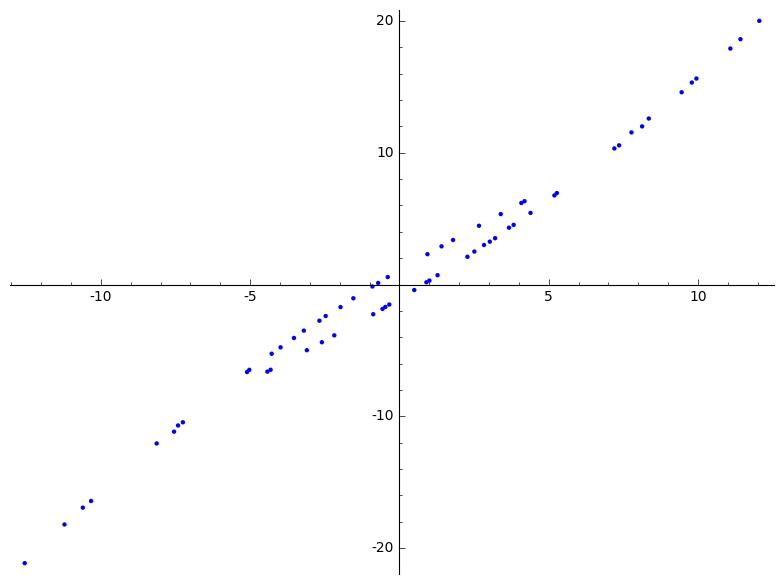}\\
\includegraphics[scale=0.25]{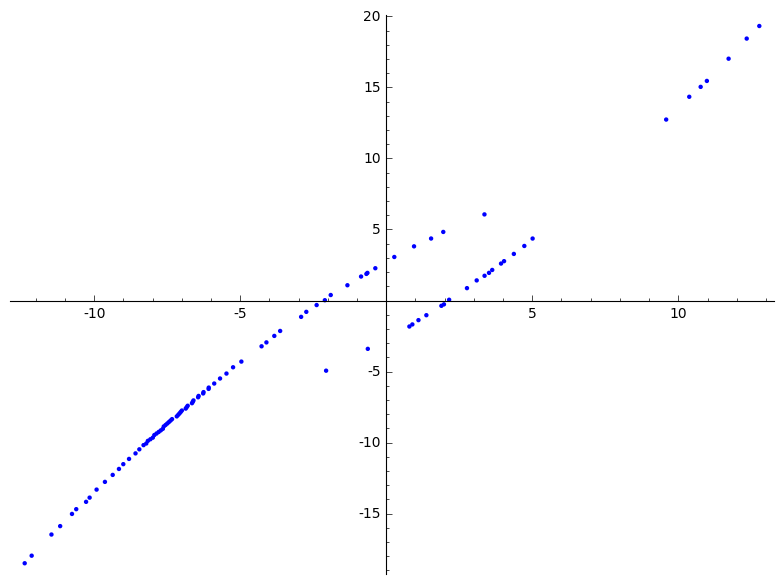} & \includegraphics[scale=0.25]{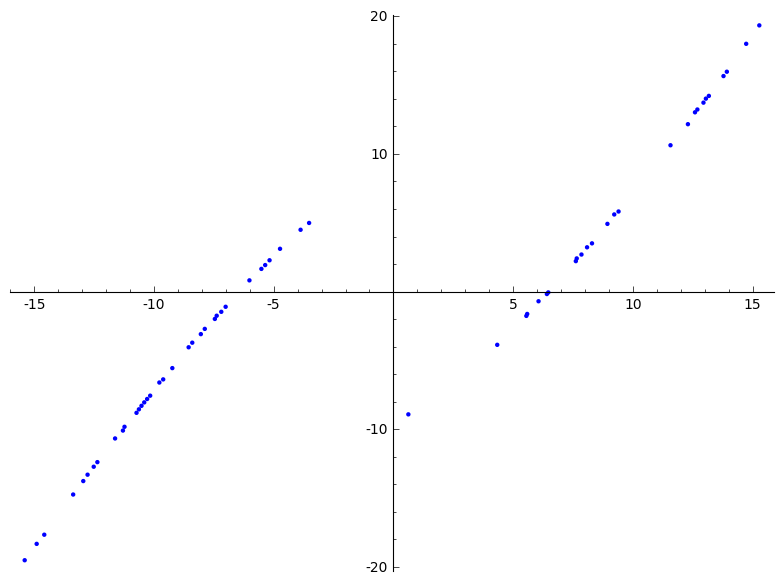}\\
\end{tabular}
\caption[AtractorSections]{Some Poincar\'e sections of Lorenz Attractor}
\end{figure}

The idea of Poincar\'e sections of Lorenz attractor will bring the next starring of our story: the french mathematician and astronomer Michel H\'enon (1931-2013). Despite his main research was in  restricted-three-body and n-body problems, he also did something quite relevant in discrete dynamics: Poincar\'e sections of Lorenz system.

\subsection{H\'enon Map}

In the paper A Two-dimensional Mapping with a Strange Attractor (\cite{Henon}) H\'enon shows that a discrete mapping of the plane holds the same properties of Lorenz attractor. In Lorenz attractor, the divergence of the flow has a constant negative value, so any volume shrinks exponentially with time, with a bounded region trapping all trajectories. H\'enon pointed out that the strange attractor seems to be (locally) the product of a two-dimensional manifold by a Cantor set.

\begin{definition}
A set $\Lambda$ is a Cantor Set if it is a closed, totally disconnected, and perfect subset of $I$. A set is totally disconnected if it contains no intervals; a set is perfect if every point in it is an accumulation point of limit point of other points in the set.
\end{definition}

So finding a simpler model of Lorenz strange attractor allows better (qualitative and quantitative) explorations of it. As a first step, Poincar\'e mapping of Lorenz attractor will be considered, instead of the whole set. This will "decrease" the problem in one less dimension. However, it still requires the numerical integration of the differential equation so a new (explicit) map $T$ will be defined. even if it is not Lorenz system, essential propertied will be held.
\\
Pomeau (1976) showed that, in Lorenz system, a volume is stretched in one direction, and simultaneously folded over itself, in the course of a revolution. Based on this H\'enon took an elongated region and applied the following transformations:

\begin{eqnarray*}
T': x'=x, y'=y+1-ax^2\\
T'': x''=bx', y''=y'\\
T''': x'''=y'', y'''=x''\\
\end{eqnarray*}

The H\'enon mapping is then defined as $T=T'''T''T'$.

\begin{definition}{Hen\'on Map}
$$T: x_{i+1}=y_{i}+1-ax_{i}^{2}, y_{i+1}=bx_{i}$$
\end{definition}

\begin{figure}
\centering
\begin{tabular}{cccc}
\includegraphics[scale=0.10]{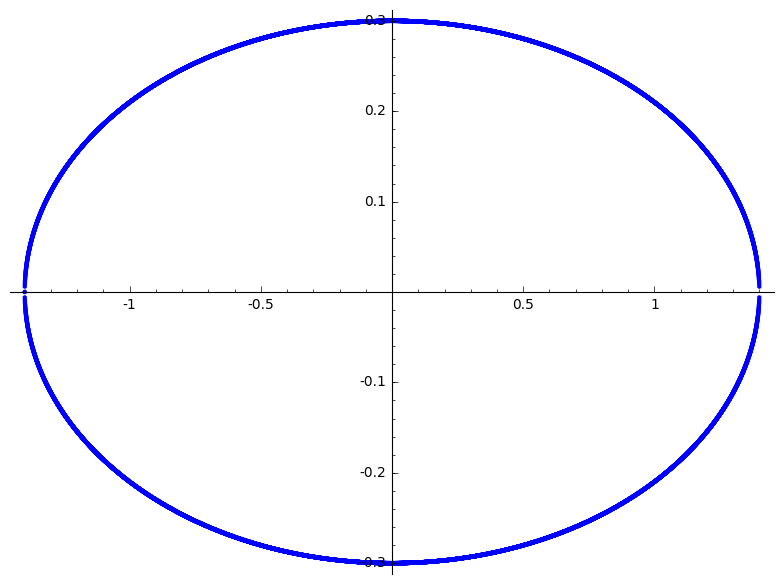} &  &  & \includegraphics[scale=0.10]{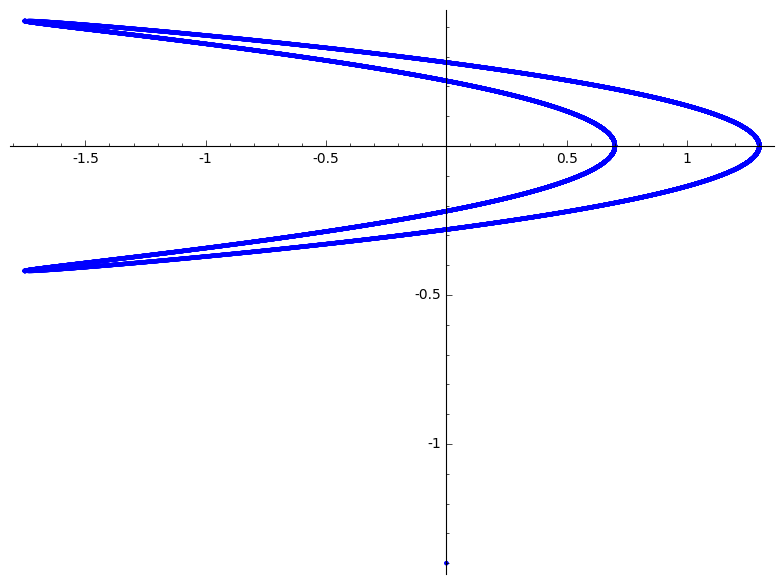} \\
 & \includegraphics[scale=0.10]{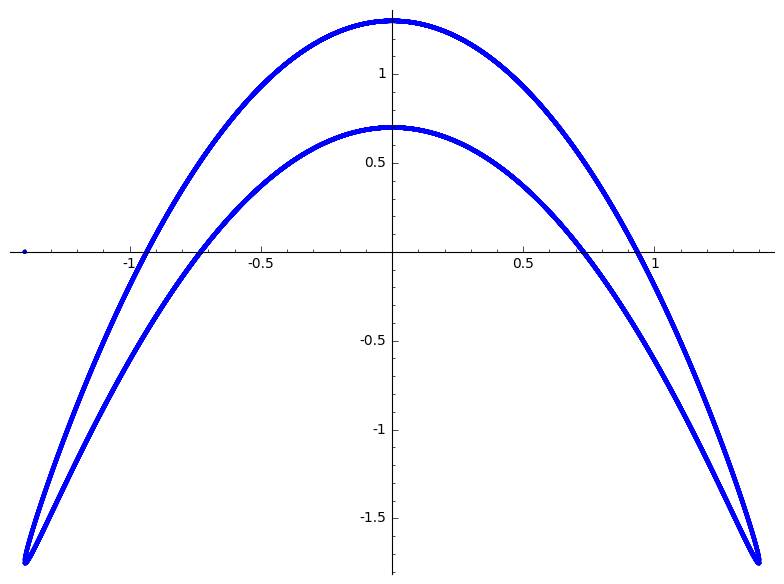} &\includegraphics[scale=0.10]{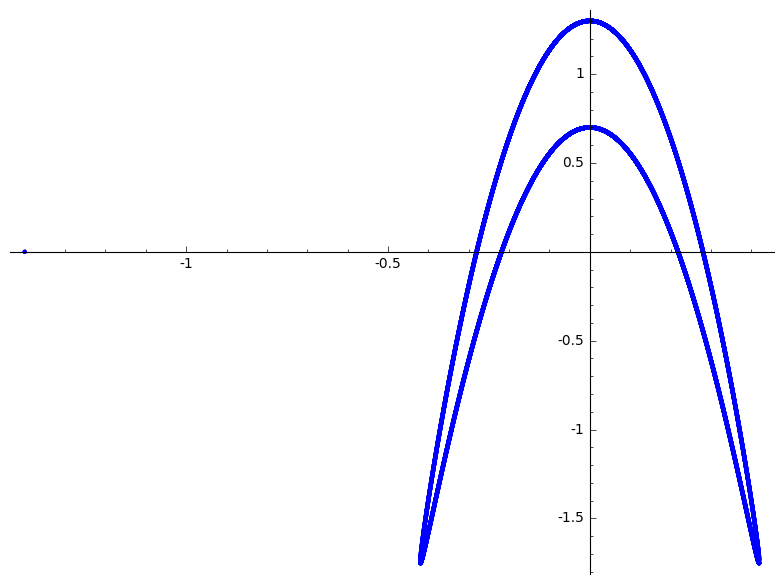} &  \\
\end{tabular}
\caption[HenonSteps]{Transformations to obtain the H\'enon map $T: x_{i+1}=y_{i}+1-ax_{i}^{2}, y_{i+1}=bx_{i}$}
\end{figure}

It is possible to state some main qualitative features of this map. The Jacobian of the H\'enon map is constant (\cite{Yoccoz}) $$\dfrac{\partial (x_{i+1},y_{i+1})}{\partial (x_{i},y_{i})}=-b$$ Besides, if $a^{2}x^{2}+b \geq 0$ the (real) eigenvalues are $$-ax \pm \sqrt{a^{2}x^{2}+b}$$ An important element for its analysis is the fact that $T$ is one to-one (\cite{LoziBook}). Thus it admits an inverse transformation $$T^{-1}: x_{i}=b^{-1}y_{i+1},  y_{i}=x_{i+1}-1+\frac{a}{b^{2}}y_{i+1}^{2}$$ Finally, $T$ has two fixed points: $x=\frac{1}{2a} \pm \sqrt{(1-b)^{2}+4a},y=bx$. The positive case is attracting if $$a_{0} = \dfrac{(1-b)^{2}}{4} < a < \dfrac{3(1-b)^{2}}{4}$$ This will be used to run numerical simulations.

In order to choose parameters, the author suggested $b=0.3$. Once this value is fixed, the following values of $a$ represent qualitative changes of behavior: $$a_{0} = \dfrac{(1-b)^{2}}{4} , a_{1} = \dfrac{3(1-b)^{2}}{4}, a_{2} \approx 1.06, a_{3} \approx 1.55$$

Taking into account that the important case of study is the strange attractor, let be $a=1.4, b=0.3$. The invariant point is approximately $(0.631354477089505 , 0.189406343126851)$, so the initial value is set as $(0.63135448 , 0.18940634)$. 

\textbf{\begin{table}
\centering
\begin{tabular}{| c | c |}
\hline
$a < a_0$ or $a > a_3$ & points escape to infinity\\
\hline
$a_0 < a < a_1$ & stable invariant point\\
\hline
$a_1 < a < a_2$ & stable invariant set of $p$ points\\
\hline
$a_2 < a < a_3$ & strange attractor\\
\hline
\end{tabular}
\caption[HenonSteps]{Behavior of $T$ respect to $a$, with $b=0.3$}
\end{table}}

\begin{figure}
\centering
\begin{tabular}{|c|c|}
\hline
\includegraphics[scale=0.2]{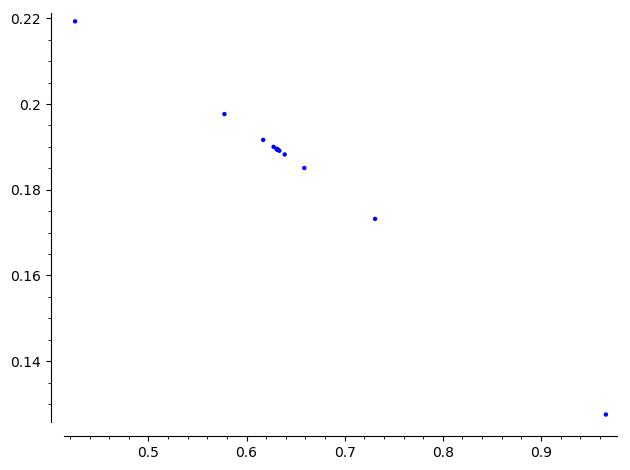} & \includegraphics[scale=0.2]{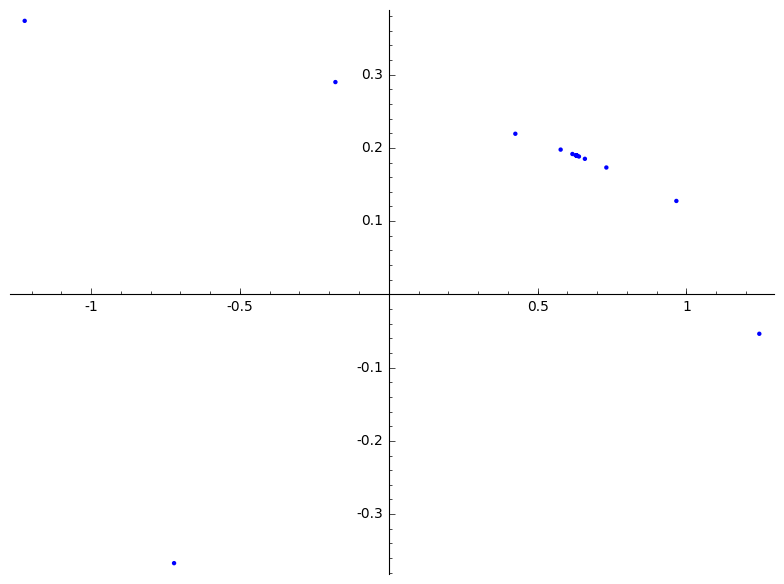}\\
\hline
\includegraphics[scale=0.2]{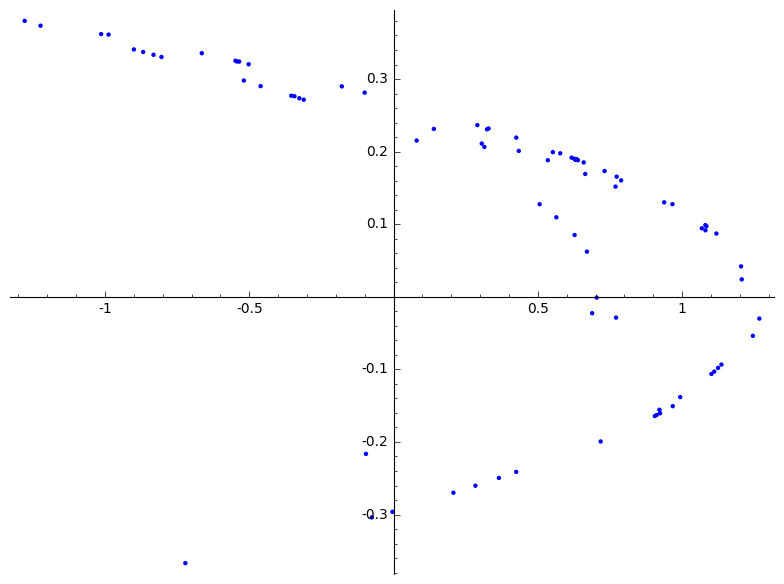} & \includegraphics[scale=0.2]{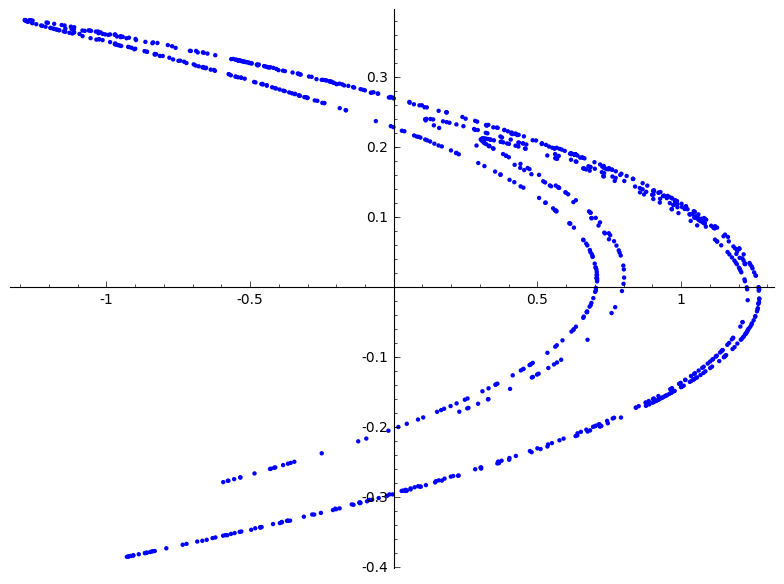}\\
\hline
\end{tabular}
\caption[HenonSteps2]{Iterates of H\'enon map $n=29,33,100,1000$}
\end{figure}

In the first numerical simulation, H\'enon used $(0,0)$ as initial value. However, after finding the invariant points, He used approximations of that point to start simulations.

\begin{figure}
\centering
\includegraphics[scale=0.45]{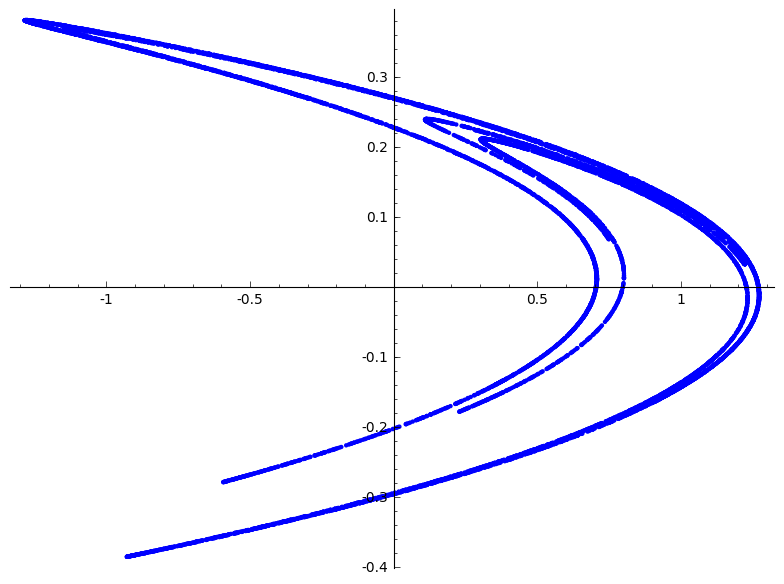}
\caption[Henon10000]{H\'enon map $n=10^4$}
\end{figure}

The eigenvalues associated to this invariant point are 
\begin{eqnarray*}
\lambda_1 = 0.15594632... \\
\lambda_2 = -1.92373886... \\
\end{eqnarray*}

The slopes given by the eigenvectors are
\begin{eqnarray*}
p_1 = 1.92373886... \\
p_2=-0.15594632... \\
\end{eqnarray*}

Due to $\lambda_2$, this point is unstable. and the lines passing through that point with slopes $p_1$ and $p_2$ are the stable and unstable manifolds respectively. The lines seem to be continuous, but after an enlargement, it is evident that each curve is made of curves, as well as its similarity with a Cantor Set

\begin{figure}
\centering
\includegraphics[scale=0.45]{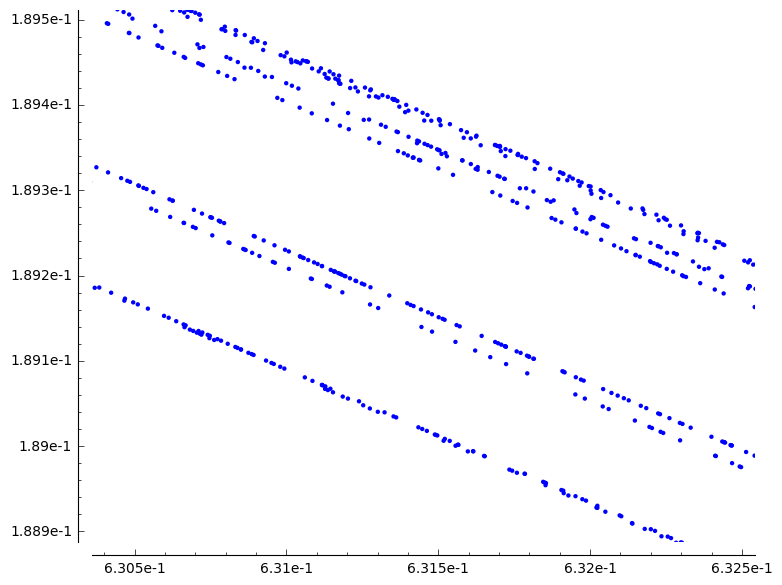}
\caption[Henon1000000]{H\'enon map $n=5 \times 10^6$, CPU time: 707.81 s}
\end{figure}

Summarizing, Hen\'on map shows a similar structure to Lorenz attractor, but its numerical exploration is simpler. Lorenz and Pomeau inferred the Cantor set, but due to the contracting ratio on Lorenz map ($7 \times 10^5$) they couldn't observe it, while $T$ has a smaller ratio ($b=0.3$). Now it is time to simplify this problem one more time, going to one dimension and face the Cantor set in a simpler way.

\section{May}

Robert M. May (1936-2020) was, in the most general sense of the word, a scientist. Despite being educated in chemistry and physics, he develop significant research in population dynamics, specially in animal populations. This work was quite relevant to develop theoretical ecology, and also states a new context to explore discrete dynamic systems: Difference equations. Based in a particular equation (known as logistic model) he introduced chaotic behaviour by double-period bifurcation. 

In Simple mathematical models with very complicated dynamics (\cite{May}), the author summarized non-linear phenomena, showing the way in which those simple models, used to model dynamics of biological populations, move from stable points to chaos. In first-order difference equations, a variable $X_{t+1}$ is related to its preceding value $X_t$. This is usually expressed as $X_{t+1} = F(X_{t})$. A particular (and relevant) case is the logistic difference equation $N_{t+1} = N_{t}(a-b N_{t})$. Taking $X=\frac{b N}{a}$, the equation turns into $X_{t+1} = a X_{t} (1+X_{t})$. This last form was used for May to develop his paper. In order to keep all iterates in $[0,1]$, the maximum value of f must be $\frac{a}{4}=1$, so $a<4$. This is interesting, because May developed bifurcation theory with this restriction over $a$, however, the case $a>4$ is as much as interesting as the previous one, due to the emerging Cantor set.

\subsection{Quadratic Map, $a<4$}

\begin{definition}
Let $f^{n}(x)$ be the $f \circ f \circ \ldots \circ f$ $n$ times.

A point $x$ is fixed for $f(x)$ if $f(x)=x$. A point $x$ is said to be periodic of period $n$ if $f^{n}(x)=x$, where the least $n$ holding this equation is named the prime period of $x$. (\cite{Devaney})
\end{definition}

Fixed point are analogous to equilibrium values in (continuous) differential equations. In this case the fixed points are:

\begin{eqnarray*}
X^{*} = a X^{*} (1-X^{*})\\
X^{*} = a X^{*} - a (X^{*})^2\\
(1-a)X^{*} = - a (X^{*})^2\\
a (X^{*})^2 + (1-a)X^{*} = 0\\
X^{*} (a X^{*} + (1-a))=0\\
X^{*}=0; X^{*}=\frac{a-1}{a}
\end{eqnarray*}

\begin{definition}
A fixed point $p$ is named hyperbolic if $\vert f'(p) \vert \neq 1$. This point will be an attracting point if $\vert f'(p) \vert < 1$, and a repelling point if $\vert f'(p) \vert > 1$
\end{definition}

Since $X'(0)=a$ and $X'(\frac{a-1}{a}) = 2 - a$. So $0$ is a repelling point for $a>1$, and $\frac{a-1}{a}$ is an attracting point for $1<a<3$. If $a>3$, the point $\frac{a-1}{a}$ is not hyperbolic, however, a pair of fixed point appear in $X^{(2)}$, it means, a period-2 point. This is known as a double-period bifurcation of the parameter $a$, and the points repelled by $\frac{a-1}{a}$ lie on the new periodic orbit.

\begin{figure}
\centering
\includegraphics[scale=0.5]{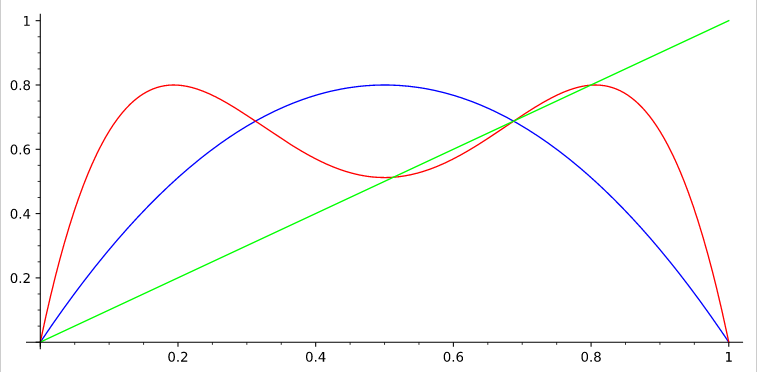}
\caption[LogisticMap32]{$X_{t}^{(2)}$ with $a=3.2$}
\end{figure}

The period-2 orbit is stable if $a < 3.449489742783179$. After this value of $a$ the period-2 orbit turns into a period-4 orbit, showing again a double-period bifurcation. This bifurcation will appear as $a$ increases, generating period-$2^n$ up to the accumulation point at $a \approx 3.57$ (\cite{May}). Beyond this point, infinite periodic points with different periods will appear. In fact, at $a \approx 3.8284$ appears a period-3 point, and after that, it is possible to find periodic points of all periods, as well as aperiodic points. This step to chaos was proofed simultaneously (and independently) by Li and Yorke in (\cite{LiYorke}) and Sharkovskii (\cite{Sharkovskii}). This last approach is one of the main referents of this work.

Still, there is one loose end: the relation between May's difference equation and H\'enon map. These two objects will meet in a common place: the Cantor Set, and this is the reason why it makes sense "to go down" from 2 to 1 dimension. in this travel from Poincar\'e to May, it is possible to see how dynamic phenomenena are quite similar in (apparently) different structures. This is the importance of Discrete Dynamics as a tool to understand continuous systems. In order to conclude this first part, the Cantor Set structure of the quadratic map will be studied.  

\subsection{Quadratic map, $a \geq 4$}

The variation of $a$ in the interval $1<a<4$ produces an amazing chain reaction of bifurcations from fixed points with stable orbits to chaos but, What if $a \geq 4$?. If  $x \in I = [0,1]$ and $1 \leq a \leq 4$, the map $X_{t+1} = a X_{t} (1+X_{t})$ sends the interval $I$ into a set $S \subset I$, keeping all iterates into $I$. On the other hand, if $a>4$, the image of $I$  is $[0,\frac{a}{4}] \supset I$, and the orbit of points in the interval $A_0 = [\frac{a-\sqrt{a(a-4)}}{2a},\frac{a+\sqrt{a(a-4)}}{2a}]$ is unstable, since $X_1 >1$, $X_2 < 0$ and $X_n \rightarrow -\infty$.

After the first iterate, only the points in $I \setminus A_0$  have an image in $I$. Let be $I_0 \cup I_1 = I \setminus A_0$ the disjoint intervals produced by subtracting $A_0$.  Analogously, the points of $I$ with image in $A_0$ in the first iterate, will be part of an unstable orbit. The inverse image of $A_0$ is a pair of intervals $A_{1,0}, A_{1,1}$ such that $A_1 = A_{1,0} \cup A_{1,1}$ and $A_{1,0} \subset I_0, A_{1,1} \subset I_1$, So after this, the points remaining in $I$ after 2 iterates are in four disjoint subsets of I, generated by $I_0 \setminus A_1$ and $I_1 \setminus A_2$, it means, taking "a middle third" of the intervals again. 

Let be $A_n = \left\lbrace  x \in I / X_{i} \in I, X_{n+1} \notin I \right\rbrace $. Points in $A_n$ eventually will be part of an unstable orbit, so the points with full dynamics in $I$ will be the points in $$\Lambda = I - \bigcup\limits_{i=0}^{\infty} A_{i}$$ This process remains the Cantor Middle-Thirds set, since each iterate takes a "middle interval"  of the remaining subintervals of $I$. 

\begin{figure}
\centering
\begin{tabular}{c c}
\includegraphics[scale=0.5]{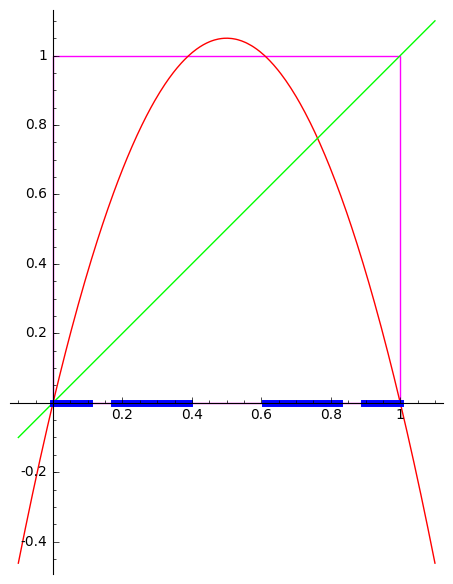} &
\includegraphics[scale=0.5]{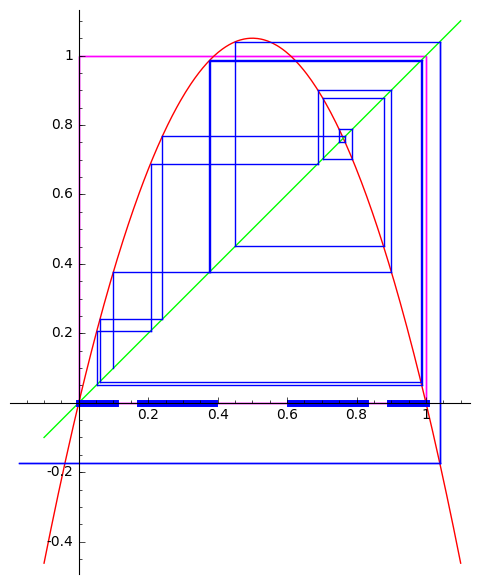} \\
\end{tabular}
\caption[CantorSet]{Iterates of $X_n$ with $a>4$}
\end{figure}

\begin{theorem}
If $a$ accomplishes that $(\forall x \in I \setminus A_0)(\vert x' \vert > 1)$ then $\Lambda$ is a Cantor set (\cite{Devaney})
\end{theorem}

Setting a proper $a$ as above, there exists $\lambda > 1$ such that $(\forall x \in \Lambda)(\vert x'_n \vert > \lambda)$, and it implies that $\vert (x_{n})' \vert > \lambda ^ n$. This is relevant to show $ \Lambda $ is totally disconnected. 

Suppose $x\neq y, [x,y] \subset \Lambda$, if $c \in [x,y]$ then $\vert (c_{n})' \vert > \lambda ^ n$. Due to archimedean property it is possible to choose $n$ so that $\lambda ^n \vert x-y \vert > 1$. According to the Mean Value Theorem $\vert y_n - x_n \vert \geq \lambda ^ n \vert y-x \vert > 1$ so either $x_n$ or $y_n$ is out of $I$, which is a contradiction, so $\Lambda$ is totally disconnected. Besides, $\Lambda$ is closed, since it is the intersection of nested closed intervals.

Finally, $\Lambda$ must be perfect. All endpoints of $A_n$ will have 0 as image, staying in $I$ from that point and on. Consider now an isolated point $c \in \Lambda$. points in a neighbourhood of $p$ should leave $I$ after certain $n*$, so these points belong to $A_{n*}$. There are two possibilities: a sequence of endpoints converges to $c$, or all points in deleted neighbourhood of $c$ leaves $I$ eventually. In the first case, the points will stay in $\Lambda$. in the second one, $X_n$ maps $c$ to 0, and its neighbourhood out of $I$ (into the negative real axis). $X_n$ has a maximum at $p$ since $(p_n)'=0$, so $p_i = \frac{1}{2}$ for some $i<n$. This implies that $p_{i+1} \notin I$, which is again a contradiction.

\section{Final Remarks}
Dynamical systems are very important in the development of different areas, see for example \cite{AAT} for applications in Quantum Mechanics and see also \cite{siam,AABD,AAD,AAS} for applications in Classical Mechanics. In particular, discrete dynamical systems play an important role due to they are an useful tool to understand continuous systems by decreasing its complexity, either by diminish dimensions or passing from continuous to discrete time. Besides, they are a rich research field, with a significant number of unsolved problems, such as genealogy of permutations, a qualitative description of forcing relationships in bifurcations, and the extension of most of the 1-dimension methods to 2-dimensional systems.

\end{document}